\documentclass{siamart1116}

\usepackage{amsmath}
\usepackage{amsfonts}
\usepackage{amssymb}
\usepackage{graphicx}

\usepackage{booktabs}

\usepackage{algorithm}
\usepackage{algpseudocode}
\let\oldReturn\Return
\renewcommand{\Return}{\State\oldReturn}
\usepackage{listings}

\newcommand{\eps}{\mbox{$\epsilon$}}

\let\vec=\mathbi%
\let\mat=\mathbf%
\let\set= \mathcal%

\usepackage{subcaption}
\usepackage{graphicx}
\usepackage{array}
\usepackage{multirow}
\newcommand{\specialcell}[2][c]{\begin{tabular}[#1]{@{}c@{}}#2\end{tabular}}
\usepackage{amsmath}
\usepackage{mleftright}
\usepackage{amssymb}
\usepackage{color}
\usepackage{tikz}
\usepackage{hhline}
\usepackage{booktabs}
\usetikzlibrary{matrix, automata}
\definecolor{RED}{rgb}{1,0,0}

\definecolor{GREEN}{rgb}{0,0,0}
\definecolor{BLUE}{rgb}{0,0.4,1}
\newcommand\given[1][]{\:#1\vert\:}

\DeclareMathOperator*{\argmin}{arg\,min}
\DeclareMathOperator*{\argmax}{arg\,max}

\usepackage{lipsum}
\usepackage{amsfonts}
\usepackage{graphicx}
\usepackage{epstopdf}
\usepackage{tensor}
\ifpdf
  \DeclareGraphicsExtensions{.eps,.pdf,.png,.jpg}
\else
  \DeclareGraphicsExtensions{.eps}
\fi

\numberwithin{theorem}{section}

\newcommand{\TheTitle}{Tensor Approximation of Advanced Metrics for Sensitivity Analysis}
\newcommand{\TheAuthors}{R. Ballester-Ripoll, E. G. Paredes, and R. Pajarola}

\headers{\TheTitle}{\TheAuthors}

\title{{\TheTitle}\thanks{ Submitted to the editors: 6th December 2017.
\funding{This work was partially supported by the University of Zurich's Forschungskredit ``Candoc'', grant number FK-16-012.}}}

\author{
  Rafael Ballester-Ripoll\thanks{Department of Informatics, University of Zurich
    (\email{rballester@ifi.uzh.ch}, \email{egparedes@ifi.uzh.ch}, \email{pajarola@ifi.uzh.ch}).}
  \and
  Enrique G. Paredes\footnotemark[2]
  \and
  Renato Pajarola\footnotemark[2].
}

\usepackage{amsopn}

\ifpdf
\hypersetup{
  pdftitle={\TheTitle},
  pdfauthor={\TheAuthors}
}
\fi

\begin{document}

\maketitle

\begin{abstract}
  Following up on the success of the analysis of variance (ANOVA) decomposition and the Sobol indices (SI) for global sensitivity analysis, various related quantities of interest have been defined in the literature including the effective and mean dimensions, the dimension distribution, and the Shapley values. Such metrics combine up to exponential numbers of SI in different ways and can be of great aid in uncertainty quantification and model interpretation tasks, but are computationally challenging. We focus on surrogate based sensitivity analysis for independently distributed variables, namely via the tensor train (TT) decomposition. This format permits flexible and scalable surrogate modeling and can efficiently extract all SI at once in a compressed TT representation of their own. Based on this, we contribute a range of novel algorithms that compute more advanced sensitivity metrics by selecting and aggregating certain subsets of SI in the tensor compressed domain. Drawing on an interpretation of the TT model in terms of deterministic finite automata, we are able to construct explicit auxiliary TT tensors that encode exactly all necessary index selection masks. Having both the SI and the masks in the TT format allows efficient computation of all aforementioned metrics, as we demonstrate in a number of example models.
\end{abstract}

\begin{keywords}
  Variance-based sensitivity analysis, surrogate modeling, tensor train decomposition, Sobol indices
\end{keywords}

\begin{AMS}
  65C20, 15A69, 49Q12
\end{AMS}

\section{Introduction}

Variance-based sensitivity analysis (SA) is a fundamental tool in many disciplines including reliability engineering, risk assessment and uncertainty quantification. It captures the behavior of simulations and systems in terms of how much of their output's variability is explained by each input (and combinations of inputs), and has received a great deal of academic and industrial interest over the last decades. These efforts have resulted in widely popular metrics such as the \emph{Sobol indices} (SI)~\cite{Sobol:90, HS:96} and an increasing number of more recent related quantities of interest (QoI). They help analysts assess which groups of variables have the strongest influence on the output's uncertainty and, for example, which ones may be frozen with the least possible impact~\cite{STCR:04}. 

A number of long-standing hurdles make such tasks challenging. First of all, directly sampling the whole domain of variables is rarely a feasible option. Usually one has either a given sparse set of fixed samples, or a simulation/experiment that can be run on demand with arbitrary parameters (for example the so-called \emph{non-intrusive modeling}, also known as \emph{black-box sampling}). Uncertainty estimations are thus often bound to have a margin of error. Second, the well-known \emph{curse of dimensionality} poses a challenge for high-parametric models. Points tend to lie \emph{far} from each other as the number of variables $N$ grows, and a rather large number of samples may be required in practice to attain a reasonable accuracy. Furthermore, the number of possible index combinations and metrics scales exponentially with $N$. Several algorithms and sampling schemes have been proposed that can partially tackle this problem; for example to estimate specific aggregated indices (e.g. the \emph{total effects}) or as exact formulas to compute analytical values from certain classes of mathematical functions. Many methods limit themselves to computing indices that are relative to single variables only. Unfortunately, such simplifications risk overlooking sizable joint interactions. Often, an effect due to a specific combination of 2 or 3 variables might be stronger and more significant than the mere knowledge that these variables are important on their own.

A more powerful strategy for SA consists in using a limited set of samples to train a regressor that acts as a \emph{surrogate model}, i.e. a routine that can estimate the true model's output for any combination of input values. This has become a standard choice for many SA tasks~\cite{QHSGVT:05, Yankov:15, LMS:17, IP:17}, especially when many samples are needed. Even though such models are typically very fast to evaluate, sampling schemes operating on them still suffer from the curse of dimensionality, and certain QoIs or queries can be highly time-consuming to compute. 
In particular, some advanced metrics are defined on and combine many or even the whole set of SI. Even though interactions involving many variables tend to be very small in practice, there is an exponentially large number of them; hence their aggregated contributions should not be ignored in general.

The present work tackles surrogate modeling based SA under the paradigm of low-rank tensor decompositions, namely the \emph{tensor train} (TT) model~\cite{Oseledets:11}. This model lends itself very well to variance-based SA as the SI can be extracted directly from its compressed representation without explicit sampling~\cite{Rai:14, DKLM:14, BPP:17}. Throughout this paper we assume that a low-rank TT tensor surrogate exists that approximately predicts the model behavior at all possible input variable combinations. This assumption holds for many families of models, including many with high orders of variable interactions, and also in the presence of categorical variables. Many algorithms have been proposed to build such TT representations, be it via fixed sets of samples~\cite{Steinlechner:15, GKK:15, GK:17}, via adaptively sampling black-box simulations and analytical functions~\cite{OT:10, SO:11, BEM:16}, or from other alternative low-rank tensor decompositions~\cite{Handschuh:15, KS:15, BPP:17}; see also~\cite{GKT:13, BPP:16}. In this work we focus on adaptive sampling by the so-called \emph{cross-approximation} technique~\cite{OT:10}.

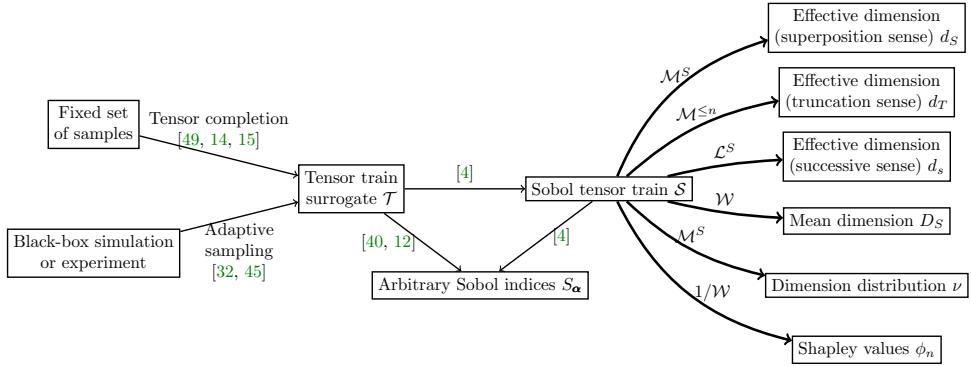
\begin{figure}[ht]
\centering
\resizebox{1\columnwidth}{!}{
\begin{tikzpicture}[scale=1.3]
\begin{scope}[every node/.style={rectangle,thick,draw}]
    \node (samples) [align=center] at (0,1) {Fixed set \\ of samples};
    \node (simulation) [align=center] at (0,-1) {Black-box simulation \\ or experiment};
    \node (tt) at (4,0) [draw, align=center]{Tensor train \\ surrogate $\mathcal{T}$};
    \node (below) at (6,-1.5) {Arbitrary Sobol indices $S_{\pmb{\alpha}}$};
    \node (stt) at (8,0) {Sobol tensor train $\mathcal{S}$};
    \node (edsuperposition) [align=center] at (12,2.5) {Effective dimension \\ (superposition sense) $d_S$};
    \node (edtruncation) [align=center] at (12,1.5) {Effective dimension \\ (truncation sense) $d_T$};
    \node (edsuccessive) [align=center] at (12,0.5) {Effective dimension \\ (successive sense) $d_s$};
    \node (meandimension) [align=center] at (12,-0.5) {Mean dimension $D_S$};
    \node (distribution) [align=center] at (12,-1.5) {Dimension distribution $\nu$};
    \node (shapley) [align=center] at (12,-2.5) {Shapley values $\phi_n$};
\end{scope}

\begin{scope}[every edge/.style={draw=black,thick}]
     \path [->] (samples) edge[above] node[align=center] {Tensor completion \\ \cite{Steinlechner:15, GKK:15, GK:17}} (tt);
     \path [->] (simulation) edge[below] node[align=center] {Adaptive \\ sampling \\ \cite{OT:10, SO:11}} (tt);
     \path [->] (tt) edge[above] node {\cite{BPP:17}} (stt);
     \path [->] (tt) edge[left] node {\cite{Rai:14, DKLM:14}} (below);
     \path [->] (stt) edge[above, bend left=25, line width=1.5pt] node {$\mathcal{M}^S$} (edsuperposition);
     \path [->] (stt) edge[right] node {\cite{BPP:17}} (below);
     \path [->] (stt) edge[above, bend left=15, line width=1.5pt] node {$\mathcal{M}^{\le n}$} (edtruncation);
     \path [->] (stt) edge[above, bend left=5, line width=1.5pt] node {$\mathcal{L}^S$} (edsuccessive);
     \path [->] (stt) edge[above, bend right=5, line width=1.5pt] node {$\mathcal{W}$} (meandimension);
     \path [->] (stt) edge[above, bend right=15, line width=1.5pt] node {$\mathcal{M}^S$} (distribution);
     \path [->] (stt) edge[right, bend right=25, line width=1.5pt] node {$1/\mathcal{W}$} (shapley);
\end{scope}
\end{tikzpicture}
}
\caption{Previous work has established the TT model as a valuable tool for sensitivity analysis. Based on this, in this paper we explicitly construct several tensors (bold edges) that select and aggregate Sobol indices in convenient ways and produce a variety of advanced QoIs (rightmost boxes). }
\label{fig:state_of_the_art}
\end{figure}

We contribute a range of procedures to compute the \emph{effective dimension} (in the \emph{superposition}, best-possible ordering \emph{truncation}, and \emph{successive} senses~\cite{CMO:97, LL:00}), the \emph{mean dimension}~\cite{CMO:97}, the full \emph{dimension distribution}~\cite{Owen:03}, and the \emph{Shapley values}~\cite{Shapley:53, Owen:14, Owen:17}. Current state-of-the-art approaches for these advanced metrics are narrow in scope and/or face important limitations:~\cite{SNS:16} resorts to randomly sampling the vast space of possible variable permutations to approximate the Shapley values;~\cite{RO:06} is able to estimate statistical moments of the dimension distribution; and~\cite{KFSM:11} approximates some effective dimensions by using bounds on related surrogate metrics, a method that is less effective for higher-order interaction models~\cite{BKS:15}. In contrast, we propose to use a highly compact data structure, the Sobol tensor train~\cite{BPP:17}. To the best of our knowledge, ours is the first framework that can obtain all these metrics in an efficient manner. Our algorithms exploit the fact that a certain class of finite automata can be compactly encoded using tensor networks, in particular including the TT format (see \cite{CB:08} for an early application of this interpretation for eigenstate energy minimization, and~\cite{Rabusseau:16} in the context of string weighted automata). Thanks to the advantageous numerical properties of the TT model, the proposed approach is flexible and scales well with the model dimensionality, also when the queried metrics are defined as aggregations of up to exponential numbers of Sobol indices. Fig.~\ref{fig:state_of_the_art} summarizes the state-of-the-art on TT-based SA and the algorithms here contributed.

The rest of the paper is organized as follows. Sec.~\ref{sec:definitions} reviews the main definitions and concepts we use, including the ANOVA decomposition, the Sobol indices, and all other sensitivity metrics considered. Sec.~\ref{sec:sobol_tensor_trains} outlines the mathematical tools that are fundamental for our algorithms: the tensor train decomposition (TT) as a framework for surrogate modeling and the Sobol tensor train. In Sec.~\ref{sec:automata} we contribute several explicitly constructed TT tensors which behave as deterministic finite automata on tensors of size $2^N$. In Sec.~\ref{sec:computation} we show how these automata can be combined with Sobol TTs to efficiently produce all sensitivity metrics listed in Sec.~\ref{sec:definitions}. Numerical results are presented in Sec.~\ref{sec:results}, and concluding remarks in Sec.~\ref{sec:conclusions}.

\section{Sensitivity Analysis: Definitions and Metrics} \label{sec:definitions}

We write multiarrays (tensors) differently depending on their dimension: italic scalars (e.g. $s$), vectors with boldface italics (e.g. $\vec{u}$), matrices with boldface capitals (e.g. $\mat{U}$), and higher-order tensors with calligraphic capitals (e.g. $\set{T}$). We point to their elements via NumPy-like notation, for example $\set{T}[:, \mathrm{:}k, :]$ are the slices $0, \dots k-1$ of a 3D tensor along its second axis. 
Tuples of model variables are denoted as $\pmb{\alpha}$, and Sobol indices as $S_{\pmb{\alpha}}$. We write complements of tuples as $-\pmb{\alpha} := \{1, \dots, N\} \setminus \pmb{\alpha}$. Often, we access tensors using variable tuples, and write them in subscripts for short: $\set{T}_{\pmb{\alpha}} \equiv \set{T}[\alpha_1, \dots, \alpha_N]$ where $\alpha_n = 1$ iff $n \in \pmb{\alpha}$, and 0 otherwise. For example, if $\set{T}$ is a 3D tensor, $\set{T}_{2, 3}$ denotes the element $\set{T}[0, 1, 1]$. The Kronecker and element-wise products for matrices and tensors are written as $\mat{U} \otimes \mat{V}$ and $\set{A} \circ \set{B}$, respectively.

\subsection{ANOVA Decomposition and Sobol Indices}

Given a function $f:\Omega \subset \mathbb{R}^N \to \mathbb{R}$ that is $L^2$-integrable w.r.t. a separable measure $F$, the \emph{ANOVA decomposition}, also known as the Sobol or Hoeffding decomposition~\cite{Hoeffding:48, Sobol:90}, partitions it in $2^N$ terms, each of which depends on a different subset of its input variables $\{1,\dots,N\}$:

\begin{equation} \label{eq:sobol_decomposition}
\begin{aligned}
f(\vec{x}) = f_{\emptyset} + f_1(x_1) + \dots + f_N(x_N) + 
\\
+ f_{1,2}(x_1, x_2) + \dots + f_{1 \dots N}(x_1, \dots, x_N) = \sum_{\pmb{\alpha} \subseteq \{1, \dots, N\}} f_{\pmb{\alpha}}(x_{\pmb{\alpha}})
\end{aligned}
\end{equation}

\noindent
This decomposition is unique provided that the subfunctions $f_{\pmb{\alpha}}$ have 0 mean for all $\pmb{\alpha} \subseteq \{1,\dots,N\}, \pmb{\alpha} \ne \emptyset$, w.r.t. the separable measure $F$:
$$\int_{\Omega_{\pmb{\alpha}}} f_{\pmb{\alpha}}(\vec{x}_{\pmb{\alpha}}) \, dF_{\pmb{\alpha}}(\vec{x}_{\pmb{\alpha}}) = 0$$


\noindent
and each component can be computed as

\begin{eqnarray}
& f_{\emptyset}(\vec{x}) = f_{\emptyset} = \mathbb{E}[f] = \int_{\Omega} f(\vec{x}) \, dF(\vec{x}) \\
& f_{\pmb{\alpha}}(\vec{x}) = f_{\pmb{\alpha}}(\vec{x}_{\pmb{\alpha}}) := \int_{\Omega_{-\pmb{\alpha}}} \left( f(\vec{x}) - \sum_{\pmb{\beta} \subsetneq \pmb{\alpha}} f_{\pmb{\beta}}(\vec{x}_{\pmb{\beta}}) \right) \, dF_{-\pmb{\alpha}}(\vec{x}_{-\pmb{\alpha}}) \nonumber
\end{eqnarray}

If $F$ is an $N$-dimensional separable joint probability density function (PDF), namely of a vector of independent random variables $\vec{x} = (x_1, \dots, x_N)$, then Eq.~\ref{eq:sobol_decomposition} defines a partition of the model's total statistical variance as the sum of variances of each subfunction:

\begin{equation} \label{eq:variance_partition}
\mathrm{Var}[f] = \sum_{\pmb{\alpha} \ne \emptyset} \mathrm{Var}[f_{\pmb{\alpha}}]
\end{equation}

The \emph{Sobol indices}, or SI for short, are the relative variance contributions of each subfunction (except $f_{\emptyset}$), normalized by the total variance:

\begin{equation} \label{eq:sobol_indices}
S_{\pmb{\alpha}} := \frac{\mathrm{Var}[f_{\pmb{\alpha}}]}{\mathrm{Var}[f]} \ \ (\mbox{for } \pmb{\alpha} \subseteq\{1, \dots, N\}, \pmb{\alpha} \ne \emptyset)
\end{equation}

All SI are nonnegative, since they are second-order moments. Also, from Eqs.~\ref{eq:variance_partition} and~\ref{eq:sobol_indices} it follows that the sum of all SI is 1. These properties make them interpretable in terms of set cardinalities, and can be used to define a set algebra. Some aggregations and combinations of SI have names of their own; most importantly the \emph{total Sobol indices} $\set{S}^T$ (sometimes also known as \emph{total effects} or \emph{first-order indices}) and the \emph{closed Sobol indices} $\set{S}^C$~\cite{Owen:13}:

\begin{equation}
S_{\pmb{\alpha}}^T := \sum_{\pmb{\beta} \cap \pmb{\alpha} \ne \emptyset} S_{\pmb{\beta}}, \ S_{\pmb{\alpha}}^C := \sum_{\pmb{\beta} \subseteq \pmb{\alpha}} S_{\pmb{\beta}}
\end{equation}

\noindent
which satisfy

\begin{equation}
\begin{gathered}
0 \le S_{\pmb{\alpha}} \le S_{\pmb{\alpha}}^C \le S_{\pmb{\alpha}}^T \le 1 \\
S^T_{\pmb{\alpha}} = 1 - S^C_{-\pmb{\alpha}} \mbox{ and } S^C_{\pmb{\alpha}} = 1 - S^T_{-\pmb{\alpha}}
\end{gathered}
\end{equation}

\subsection{Effective Dimension}

More advanced sensitivity metrics take into account the size of the variable tuple $\pmb{\alpha}$. One example is the \emph{effective dimension}, which has been defined in three ways at least:

\begin{itemize}
	\item \emph{Superposition sense}~\cite{CMO:97}:
	$$
	d_S := \argmin_k \left\{ k \given[\Big] \sum_{\pmb{\alpha} \given[\big] |\pmb{\alpha}| \le k} S_{\pmb{\alpha}} \ge 1-\eps\right\}
	$$
	
	\item \emph{Truncation sense}~\cite{CMO:97}:
	$$
	d_T := \argmin_k \left\{ k \given[\Big] \\ S^C_{\{1 \dots k\}} \ge 1-\eps \right\}
	$$
	
	\item \emph{Successive sense}~\cite{LL:00}:
	$$
	d_s := \argmin_k \left\{ k \given[\Big] \sum_{\pmb{\alpha} \given[\big] \mathrm{len}(\pmb{\alpha}) \le k} S_{\pmb{\alpha}} \ge 1-\eps \right\}
	$$
\end{itemize}

In the above equation, $\mathrm{len}(\pmb{\alpha}) := \argmax_n \{n \in \pmb{\alpha}\} - \argmin_n \{n \in \pmb{\alpha}\} + 1$ and $\eps$ is a small tolerance for unexplained effects (say, 0.05 or 0.01). The superposition sense $d_S$ measures the minimal order of interactions needed to capture most of the model variability. In other words, it means that the model $f$ is roughly a sum of $d_S$-dimensional subfunctions; interactions involving higher numbers of variables may be safely discarded. On the other hand, $d_T$ is the number of \emph{leading} variables needed to capture a $1-\eps$ fraction of the variance. If we allow reordering the variables, $d_T$ is the minimal integer such that there exists \emph{some} tuple $\pmb{\alpha}$ with $|\pmb{\alpha}| = d_T$ such that $S^C_{\pmb{\alpha}} \ge 1-\eps$. Last, $d_s$ is informative when all variables have an inherent ordering, for example in a time series; it means that the model consists of subfunctions that depend on \emph{neighboring} variables only.

\subsection{Mean Dimension}

The \emph{mean dimension} in the superposition sense~\cite{CMO:97} is the expected value of $|\pmb{\alpha}|$, if one were to select $\pmb{\alpha}$ with probability proportional to its variance:

\begin{equation}
D_S := \sum_{\pmb{\alpha}} S_{\pmb{\alpha}} \cdot |\pmb{\alpha}|
\end{equation}

With 3 variables, for example, $D_S = S_1 + ... + 2 \cdot S_{1,2} + ... + 3 \cdot S_{1,2,3}$. This metric is a non-integer number, unlike the effective dimension, and measures the average complexity or dimensionality of a model. A result by Liu and Owen~\cite{RO:06} states that the mean dimension equals the sum of all first-order total SI: $D_S = \sum_{n=1}^{N} S^T_n$.

\subsection{Dimension Distribution}

Denoted as $\nu$, it was defined by Owen~\cite{Owen:03} as the probability mass function of the random variable $|\pmb{\alpha}|$, if one were to select $\pmb{\alpha}$ as described before. It is a discrete variable over the domain $\{1, \dots, N\}$, and each value $1 \le n \le N$ has probability

\begin{equation}
\nu(n) = \sum_{\pmb{\alpha} \given[\big] |\pmb{\alpha}| = n} S_{\pmb{\alpha}}
\end{equation}

Its expected value is the mean dimension $D_S$. Also, knowing the dimension distribution allows a direct computation of the effective dimension in the superposition sense: if we write $\bar{\nu} := \mathrm{CDF}(\nu)$, then $d_S = \lceil \bar{\nu}^{-1}(1-\eps) \rceil$. 

\subsection{Shapley Values}

They originated in game theory~\cite{Shapley:53} to determine the just retributions that each individual player $1, ..., N$ should receive from a set of potential coalitional tasks:

\begin{equation}
\phi_n := \sum_{\pmb{\alpha} \subseteq -\{n\}} \frac{|\pmb{\alpha}|! (N - |\pmb{\alpha}| - 1)!}{N!} \cdot (C_{\pmb{\alpha} \cup \{n\}} - C_{\pmb{\alpha}})
\end{equation}

Here $C_{\pmb{\alpha}}$ represents the productivity or goodness of each coalition $\pmb{\alpha}$. Shapley values have recently been reinterpreted as variance contributions in the context of SA, whereby a connection with Sobol indices was established~\cite{Owen:14}: If the productivity of each subset of variables $\pmb{\alpha}$ is taken to be $C_{\pmb{\alpha}} := S^C_{\pmb{\alpha}}$, then the $n$-th Shapley value can be computed with the simpler formula

\begin{equation}
\phi_n = \sum_{\pmb{\alpha} | n \in \pmb{\alpha}} \frac{S_{\pmb{\alpha}}}{|\pmb{\alpha}|}
\end{equation}

For example, $\phi_1 = S_1 + S_{1,2}/2 + S_{1,3}/2 + S_{1,2,3}/3$. This interpretation has been further investigated in various settings~\cite{SNS:16, IP:17, Owen:17} and found to be a good compromise between the finely-granular standard indices $S$ and the more coarse total indices $S^T$: $S_n \le \phi_n \le S^T_n \,\, \forall n = 1, \dots, N$. It has been extended for the case of dependent input variables as well. The Shapley values $\phi$ also have the desirable property that their sum equals 1.0~\cite{Owen:14}, unlike for example the total effects.

\section{Sobol Tensor Trains} \label{sec:sobol_tensor_trains}

\begin{figure}[ht]
\centering
\includegraphics[width=0.9\columnwidth]{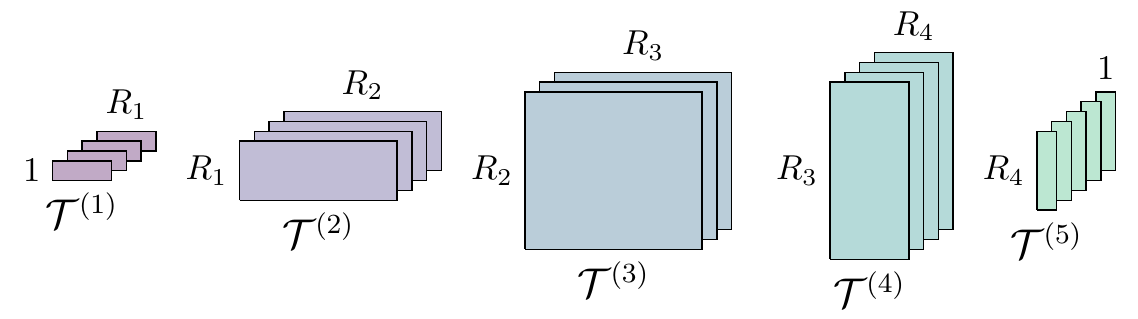}
\caption{A 5D tensor train approximating a tensor of spatial sizes $4 \times 4 \times 3 \times 4 \times 5$. Graphically, these sizes are the number of core slices (matrices) along the depth dimension, while the TT ranks $(1, R_1, R_2, R_3, R_4, 1)$ are distributed horizontally and vertically as matrix sizes. The ranks are usually larger around the central dimensions.}
\label{fig:tt}
\end{figure}

As argued in the introduction, surrogate models are a widespread approach to sensitivity analysis. In recent years, a new family of surrogates has been proposed that relies on the tensor train model (TT), a tensor decomposition proposed by Oseledets~\cite{Oseledets:11} whose size grows linearly w.r.t. the number of dimensions $N$. It is also known as \emph{linear tensor network}~\cite{Handschuh:15, CLOPZM:16} since it assigns each physical dimension to one 3D \emph{core tensor} (see Fig.~\ref{fig:tt}). To decompose a model $f$ into a tensor $\set{T}$ of size $I_1 \times \dots \times I_N$, we first discretize each axis of the domain $\Omega = \Omega_1 \times \dots \times \Omega_N$ so that $f(\vec{x}) \approx \set{T}[\vec{i}]$ with $0 \le i_n < I_n$ for $n = 1, \dots, N$. Element-wise, each entry of $\set{T}$ is a product of matrices:

\begin{equation}
\begin{split}
\set{T}[\vec{i}] = \set{T}^{(1)}[0, i_1, :] \cdot \set{T}^{(2)}[:, i_2, :] \cdot \ldots \cdot \set{T}^{(N-1)}[:, i_{N-1}, :] \cdot \set{T}^{(N)}[:, i_N, 0] = \\
\sum_{\vec{r} = \vec{0}}^{\vec{R} - \vec{1}} \set{T}^{(1)}[0, i_1, r_1] \set{T}^{(2)}[r_1, i_2, r_2] \ldots \set{T}^{(N-1)}[r_{N-2}, i_{N-1}, r_{N-1}] \set{T}^{(N)}[r_{N-1}, i_N, 0] \nonumber
\end{split}
\end{equation}

Every core $\set{T}^{(n)}$ is thus a collection of matrices that are stacked along its second dimension (corresponding to the slices in Fig.~\ref{fig:tt}). The matrix sizes are known as \emph{TT ranks} and capture the complexity of the compressed tensor. We sometimes write $\set{T} = [[\set{T}^{(1)}, \dots, \set{T}^{(N)}]]$ to denote a TT decomposition in terms of its cores.

The TT format is a powerful generalization of the SVD that can compactly encode the multidimensional structure of a wide family of functions and models. It has thus become a successful tool for high-dimensional interpolation and integration in physics and chemistry, and even in low dimensions via the so-called \emph{tensorization} process~\cite{OT:11} (which introduces artificial dimensions in a vector via reshaping operations). Given a black-box routine, an accurate TT representation can often be built using an adaptive sampling scheme over a structured set of samples, namely \emph{cross-approximation}~\cite{OT:10, SO:11}. Besides surrogate modeling~\cite{VDSL:14, BPP:16}, this model has been also used for sensitivity analysis as well~\cite{DKLM:14, Rai:14, ZYOKD:15, BEM:16, BPP:17}. For a more in-depth review on TT model building techniques from either fixed sets of samples or black-box settings, we also refer the reader to the survey~\cite{GKT:13}.

\begin{figure}[ht]
\centering
\includegraphics[width=0.9\columnwidth]{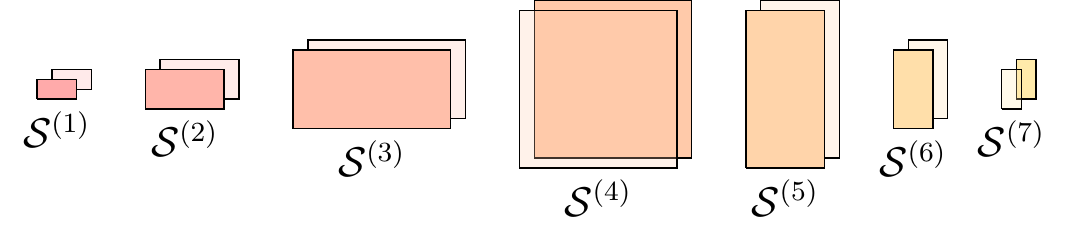}
\caption{A 7-variable model yields a Sobol TT $\set{S}$ of size $2^7$~\cite{BPP:17}. As an example, multiplying the 7 highlighted slices yields the element $\set{S}[0,0,0,1,0,0,1] = \set{S}_{4,7} \approx S_{4,7}$.}
\label{fig:sobol_tt}
\end{figure}

Ballester-Ripoll et al.~\cite{BPP:17} introduced the \emph{Sobol tensor train}, a compressed TT tensor that can be extracted from any $N$-variable TT surrogate model and approximately represents all its $2^N - 1$ Sobol indices. Denoted by $\set{S}$, the index for any tuple $\pmb{\alpha}$ is decompressed as

\begin{equation}
S_{\pmb{\alpha}} \approx \set{S}_{\pmb{\alpha}} = \set{S}[\alpha_1, \dots, \alpha_N] = \set{S}^{(1)}[:, \alpha_1, :] \cdot \ldots \cdot \set{S}^{(N)}[:, \alpha_N, :],
\end{equation}
with $\set{S}_{\emptyset} = \set{S}[0, ..., 0] = 0$. This tensor has size $2^N$ (i.e. each $\alpha_n$ can only take values in $\{0, 1\}$) and therefore each core has just 2 slices, as illustrated in Fig.~\ref{fig:sobol_tt} for a 7D Sobol TT indexing example.

Throughout this paper we will work with TT surrogate models and their Sobol TTs obtained via the method described in~\cite{BPP:17}. The Sobol TT is, to our advantage, a highly compact representation for the complete set of SIs. Furthermore, many aggregation and manipulation operations can be performed in the TT compressed domain at a cost that depends only linearly on the number of variables $N$. These include linear combinations, differentiation/integration, element-wise functions, global optimization, etc.

\section{Tensor Train Masks} \label{sec:automata}

The sensitivity metrics we overviewed in Sec.~\ref{sec:definitions} rely on selecting and weighting various SIs according to their tuple order $|\pmb{\alpha}|$. We observe that these orders can be appropriately and compactly accounted for in the TT format, thanks to a connection between tensor networks and finite automata~\cite{CB:08, Rabusseau:16}. Next we contribute the explicit construction of a number of weighted tensor masks and automata followed by their application for the direct computation of advanced metrics as further detailed in Sec.~\ref{sec:computation}. All proposed TT tensors are $N$-dimensional, i.e. use $N$ cores.

\subsection{Hamming Weight Tensor} \label{sec:hamming_weight}

We define first the \emph{Hamming weight tensor} $\set{W}$, which stores at each position $\pmb{\alpha} \in \{0, 1\}^N$ the number of `1' elements in $\pmb{\alpha}$:

\begin{equation}
\set{W}_{\pmb{\alpha}} := |\pmb{\alpha}| = \sum_{n=1}^N \alpha_n
\end{equation}

\noindent
Note that $\set{W}$ can be written as the sum of $N$ separable terms:

\begin{equation} \label{eq:bitcount}
\set{W} = \sum_{n=1}^N \overbrace{\begin{pmatrix}1 \\ 1\end{pmatrix} \otimes ... \otimes \begin{pmatrix}1 \\ 1\end{pmatrix}}^{n-1\mbox{ terms}} \otimes \begin{pmatrix}0 \\ 1\end{pmatrix} \otimes \overbrace{\begin{pmatrix}1 \\ 1\end{pmatrix} \otimes ... \otimes \begin{pmatrix}1 \\ 1\end{pmatrix}}^{N-n\mbox{ terms}}
\end{equation}

In the TT format, $\set{W}$'s rank is just 2. To see why, let us first write the Hamming weight function explicitly as

\begin{eqnarray}
\set{W}[\alpha_1, ..., \alpha_N] = \alpha_1 + ... + \alpha_N \\
\alpha_n \in \{0, 1\}, 1 \le n \le N \nonumber
\end{eqnarray}

\noindent
We now use the identity~\cite{Oseledets:13}

\begin{equation}
\label{eq:tt_identity}
\sum_{n=1}^N \alpha_n = 
\begin{pmatrix} \alpha_1 & 1\end{pmatrix}
\cdot
\begin{pmatrix} 1 & 0 \\ \alpha_2 & 1 \end{pmatrix}
\cdots
\begin{pmatrix} 1 & 0 \\ \alpha_{N-1} & 1 \end{pmatrix}
\cdot
\begin{pmatrix} 1 \\ \alpha_N\end{pmatrix}
\end{equation}

\begin{figure}[ht]
\centering
\includegraphics[width=0.5\columnwidth]{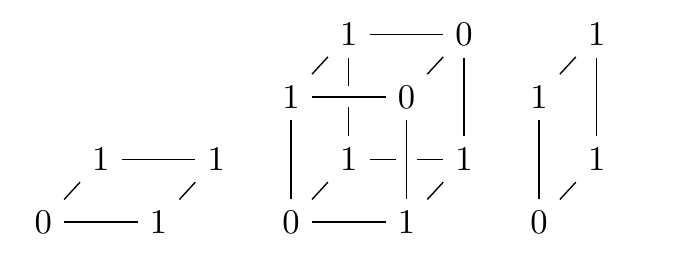}
\caption{The Hamming weight TT $\set{W}$ for $N = 3$. At each position $\pmb{\alpha} \in \{0, 1\}^N$ it records the bit sum $\set{W}_{\pmb{\alpha}} = |\pmb{\alpha}|$. It is compressed with $N$ cores (rank 2) using $8 N - 8$ elements in total.}
\label{fig:hamming_weight}
\end{figure}

If all $\{\alpha_n\}_n$ are simple scalars, Eq.~\ref{eq:tt_identity} is a ``dummy'' TT representation of a tensor of size $1^N$, with rank 2 everywhere. If every $\alpha_n$ can take values in $\{0, 1\}$, we can adapt Eq.~\ref{eq:tt_identity} by substituting each $n$-th matrix for a 3D core of size $2 \times 2 \times 2$; see Fig.~\ref{fig:hamming_weight} for an example illustration of a 3-variable model.

\subsection{Hamming Mask Tensor} \label{sec:hamming_mask}

The \emph{Hamming mask tensor} of order $\le n$ contains a $1$ at each entry $\pmb{\alpha}$ if and only if $|\pmb{\alpha}| \le n$, and a $0$ otherwise. We construct it as a deterministic finite automaton (DFA) that reads exactly $N$ symbols from the input alphabet $\{$`0', `1'$\}$ and has $n+2$ possible states $\{s_0, s_1, \dots, s_n, R\}$. There is an accepting state $s_0, \dots, s_n$ per each value of $|\pmb{\alpha}|$ for which $|\pmb{\alpha}| \le n$ and one extra rejecting state, $R$, for any other value $|\pmb{\alpha}| > n$.

\begin{figure}[ht]
\centering
\includegraphics[width=0.9\columnwidth]{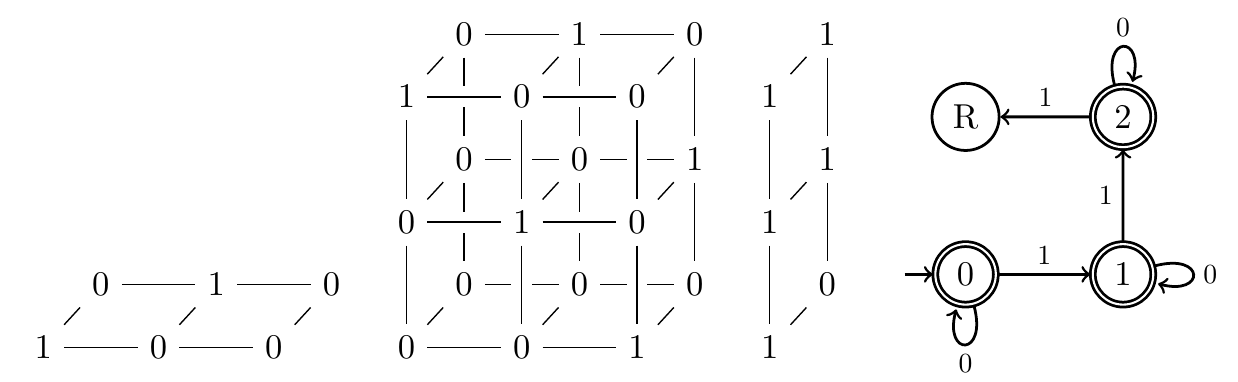}
\caption{The Hamming mask TT $\set{M}^{\le n}$ for $N = 3$ and $n = 2$ (left) and an equivalent DFA that reads 3 symbols (right, `R' stands for the rejecting state). It contains a $1$ at every position $\pmb{\alpha} \in \{0, 1\}^N$ such that $|\pmb{\alpha}| = n$, and a $0$ everywhere else. It has rank $n+1$ and uses $2 (n+1)^2 (N-2) + 4(n+1)$ elements.}
\label{fig:hamming_mask}
\end{figure}

To construct a tensor equivalent to this automaton we need to mark one of the $n+2$ possible states as the active one at each processing step. We represent the active state with a vector of $n+1$ binary elements by using \emph{dummy encoding}, with the rejecting state $R$ encoded as the all-zeros vector.
Each core updates the current state by multiplying the state vector generated in the previous step with the slice indexed by the corresponding input symbol. Thus, each slice encodes the state transition matrix for one of the input symbols (`0' or `1'). Using these conventions we manage to translate in a comprehensible manner our \emph{Hamming mask} automaton into a \emph{Hamming mask tensor} (see also Fig.~\ref{fig:hamming_mask}):

\begin{itemize}
	\item The state of the automaton after processing one input symbol can only be $s_0$ if the processed symbol was `0' (at $[0, 0, 0]$), or $s_1$ if it was `1' (at $[0, 1, 1]$). Therefore, the first core is all zeros in each slice except at those state positions.
	\item The cores $2, \dots, N-1$ have as first slice (input `0') an identity matrix so as to leave the state unchanged, while the second slice (input `1') is the identity shifted to the right since it is encoding the transition from $s_j$ to $s_{j+1}$ by shifting any `1' in the input state vector (bit at position $j$) one position to the right (to $j+1$). An intuitive explanation for these transition matrices is shown in the DFA state diagram of Fig.~\ref{fig:hamming_mask}, where transitions labeled with the `0'  input symbol always point to the current state and the ones labeled `1' point to the next state.
	\item The last core collapses the state vector so that the result is `1' if and only if the final state is one of the accepting states $\{s_0, s_1, ..., s_{n}\}$ after processing the last input symbol, that is, if $n$ or fewer `1' bits were read. That means that for the last symbol being a `0' any `1' in the state vector indicates an acceptable state $\{s_0, ..., s_{n}\}$, but for the last symbol being a `1' the state vector cannot have a `1' in the $s_{n}$-th position. Hence the first slice is all ones for checking the current state being within $\{s_0, ..., s_{n}\}$, and the second slice is all ones but at the last position to indicate the transition from $s_{n}$ to $R$.
\end{itemize}

This pattern for building a \emph{Hamming mask tensor} can be generalized to any length $N$ and any value $n$ by changing the number of cores to $N$ and their rows and columns to $n+1$ accordingly, while always keeping 2 slices.

\subsection{Hamming State Tensor} \label{sec:hamming_state}

The last core of the order $\le n$ mask tensor just described collapses the state vector into one scalar (namely, 0 or 1) and depends on $n$. However, it is also useful to output the Hamming weight using an explicit one-hot encoded result vector. We define the \emph{Hamming state tensor} $\set{M}^S$, which maps every tuple $\pmb{\alpha}$ to a vector of $N+1$ elements $\vec{u}$ as:

\begin{equation}
u_i = 
\begin{cases}
1 & \mbox{if } i = |\pmb{\alpha}| \\
0 & \mbox{otherwise}
\end{cases}
\end{equation}

To build it we proceed similarly to $M^{\le N}$, and change the last core for a 3D one, namely the same as the central ones. So $\set{M}^S$ has ranks $1, N+1, \dots, N+1, N+1$, i.e. it is not a standard TT in the sense that its last rank is not 1.

\subsection{Length Tensors} \label{sec:hamming_length}

The effective dimension in the successive sense motivates us to build tensors that are sensitive to tuple length $|\pmb{\alpha}|$, i.e. the distance between the first and the last variables present in the tuple. We define first the \emph{length mask tensor} $\set{L}^{\le n}$, which filters tuples of variables whose length exceeds a threshold (i.e are too spread out). Element-wise,

\begin{equation}
\set{L}^{\le n}_{\pmb{\alpha}} := 
\begin{cases}
1 & \mbox{if } \mathrm{len}(\pmb{\alpha}) \le n \\
0 & \mbox{otherwise}
\end{cases}
\end{equation}

\begin{figure}
\includegraphics[width=0.9\columnwidth]{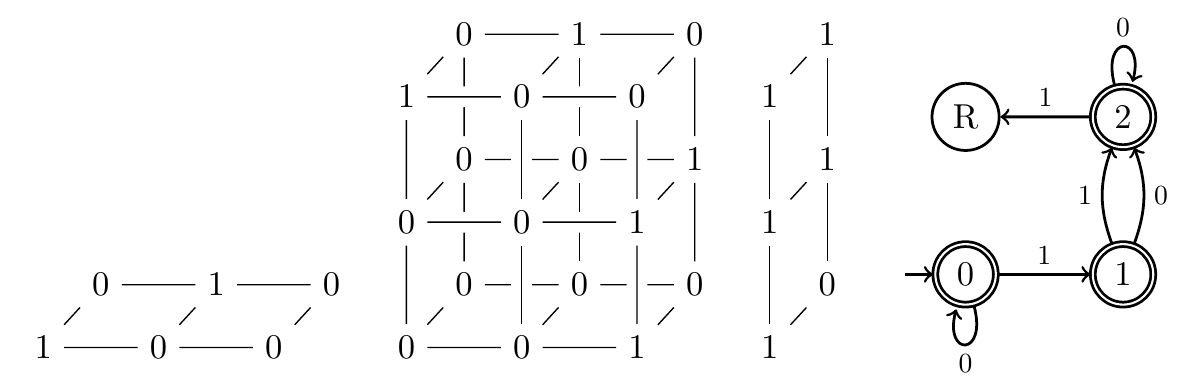}
\caption{The length mask tensor $\set{L}^{\le n}$ for $N=3$ and $n=2$ (left) and an equivalent DFA that reads 3 symbols (right). It contains a 1 at every position $\pmb{\alpha}$ where $\mathrm{len}(\pmb{\alpha}) \le n$, and 0 otherwise.}
\label{fig:lness_le_mask}
\end{figure}

We need again a state vector with $n+1$ elements to encode a state set with $n+2$ states $\{s_0, s_1, \dots, s_n, R\}$ using dummy encoding. As before, the rejecting state $R$ is represented by an all zeros state vector. Any state transition will lead to $R$ if a tuple length longer than $n$ is detected. An example of this tensor encoded automaton ($n = 2, N = 3$) is shown in Fig.~\ref{fig:lness_le_mask}:
\begin{itemize}
	\item The first core initializes the state vector in the same way as above in Sec.~\ref{sec:hamming_mask}.
	\item The cores $2, \dots, N-1$ follow a regular pattern: the second slice is the identity shifted to the right and therefore encodes the $s_j$ to $s_{j+1}$ transition for any input state. The first slice represents the transitions for the input symbol `0', which in this case are more complex:
	\begin{enumerate}
		\item When the current state is $s_0$ (a `1' in the first position of the state vector), the output state should remain the same and thus a `1' is placed at $[0, 0, 0]$ in the first row. This means that no input symbol `1' has been yet found and thus the length counter remains at $0$.
		\item When the current state is $s_1, \dots, s_{n-1}$ (a `1' in the second position of the state vector in the example), the output state should be shifted to the next one, exactly in the same way as for an input `1', thus same row as for the second slice.
		\item When the current state is $s_{n}$ (a `1' in the third position of the state vector in the example), the output state also remains $s_n$. This encodes the case in which the distance to the first input symbol `1' is already larger than $n$ and therefore the condition will be violated as soon as an input symbol `1' appears. However, if the remaining input symbols are all `0' it would mean that the actual distance to the last input symbol `1' was shorter or equal to $n$ and thus the final output would be `1'.
	\end{enumerate}
	\item The last core collapses the state vector to '1' if the final state belongs to the set of accepted states in the same way as above in Sec.~\ref{sec:hamming_mask}.
\end{itemize}

Based on this tensor and analogously to the Hamming state tensor (Sec.~\ref{sec:hamming_state}), we define now the \emph{length state tensor} $\set{L}^S$ which maps every $\pmb{\alpha}$ to a vector of $N+1$ elements that stores the value $\mathrm{len}(\pmb{\alpha})$ using one-hot encoding.

Tab.~\ref{tab:tensor_examples} shows a few examples of the values taken by the main tensors we have introduced in this section.

\begin{table*}[ht]
	\centering
	\caption{Several tensors for the case $N = 5$ evaluated at three example tuples.}
	\begin{tabular}{c c c c c c c}
		\textbf{Tuple} $\pmb{\alpha}$ & \specialcell{\textbf{Binary} \\ \textbf{form}} & $\set{W}_{\pmb{\alpha}}$ & $\set{M}^{\le 3}_{\pmb{\alpha}}$ & $\set{M}^S_{\pmb{\alpha}}$ & $\set{L}^{\le 3}_{\pmb{\alpha}}$ & $\set{L}^S_{\pmb{\alpha}}$ \\
		\hhline{=======}
		$\{5\}$ & $[0, 0, 0, 0, 1]$ & $1$ & $1$ & $[1, 0, 0, 0, 0]$ & $1$ & $[1, 0, 0, 0, 0]$ \\
		\hline
		$\{1, 5\}$ & $[1, 0, 0, 0, 1]$ & $2$ & $1$ & $[0, 1, 0, 0, 0]$ & $0$ & $[0, 0, 0, 0, 1]$ \\
		\hline
		$-\{3\} = \{1, 2, 4, 5\}$ & $[1, 1, 0, 1, 1]$ & $4$ & $0$ & $[0, 0, 0, 1, 0]$ & $0$ & $[0, 0, 0, 0, 1]$ \\
		\hline		
	\end{tabular}
	\label{tab:tensor_examples}
\end{table*}

\section{Computing Sensitivity Metrics} \label{sec:computation}


%

%
%

In this section we show how the proposed tensors can be used to extract various sensitivity analysis metrics from any Sobol TT via efficient operations in the TT format (most importantly, tensor dot products and global optimization).

\subsection{Mean Dimension} \label{sec:mean_dimension}

Let us consider the formula for the mean dimension $D_S := \sum_{\pmb{\alpha}} S_{\pmb{\alpha}} \cdot |\pmb{\alpha}|$, and let our Sobol tensor train $\set{S}$ contain all Sobol indices $S$. We observe that $D_S$ equals the tensor dot product $<\set{S}, \set{W}>$. We compute it by successively contracting the TT cores of $\set{S}$ and $\set{W}$ together (Alg.~\ref{alg:mean_dimension}, see also~\cite{Oseledets:11}) for a total asymptotic cost of $O(N R^2)$ operations, where $R$ is $\set{S}$'s rank.

\begin{algorithm}[ht]
\begin{algorithmic}[1]
\State Assemble $\set{W}$ with ranks $Q_0, Q_1, \dots , Q_{N-1}, Q_N = 1, 2, \dots, 2, 1$ (Sec.~\ref{sec:hamming_weight})
\State $\set{C} := \mathrm{ones}(1 \times 1)$ \Comment{$\set{C}$ has size $1 \times 1 = R_0 \times Q_0$}
\For{$n = 1,\dots,N$}
	\State $\set{C}_{ijk} := \set{C}_{lk} \cdot \set{S}^{l\;(n)}_{\;ji}$ \Comment{$\set{C}$ has now size $R_n \times 2 \times Q_{n-1}$}
	\State $\set{C}_{ij} := \set{C}_{imk} \cdot {\set{W}_{k\;\;j}^{\;m\;(n)}}$ \Comment{$\set{C}$ has now size $R_n \times Q_n$}
\EndFor \Comment{$\set{C}$ has now size $R_N \times Q_N = 1 \times 1$}
\Return $D_S = \mathrm{squeeze}(\set{C})$
\end{algorithmic}
\caption{Given a Sobol TT $\set{S}$, compute the mean dimension $D_S$. We use Einstein notation, i.e. tensors are contracted along the indices that appear both as subscripts and superscripts.}
\label{alg:mean_dimension}
\end{algorithm}

\subsection{Dimension Distribution}

%

\begin{figure}[ht]
\centering
\resizebox{1\columnwidth}{!}{
\includegraphics[width=0.9\columnwidth]{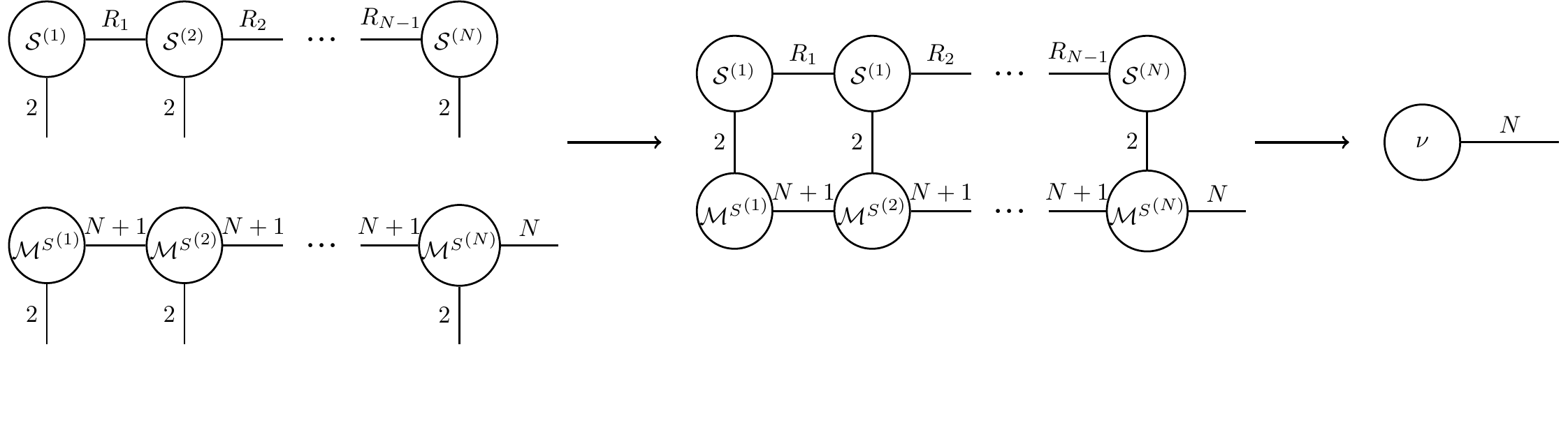}
}
\caption{We compute the dimension distribution vector $\nu$, which has $N$ elements, by contracting together the Sobol TT $\set{S}$ (ranks $1, R_1, \dots, R_{N-1}, 1$) with our proposed mask state tensor $\set{M}^S = [[\set{M}^{S(1)}, \dots, \set{M}^{S(N)}]]$ (ranks $1, N+1, \dots, N+1, 1$) along their spatial and rank dimensions. All TT cores are collapsed together, along Alg.~\ref{alg:mean_dimension}; this can be accomplished in $O(N^2 R^2) + O(N^3 R)$
operations. Only the $N$-sized edge from $\set{M}^S$ remains, and after the contraction it gathers all Sobol contributions according to their order as expected.}
\label{fig:contraction}
\end{figure}

We can extract the complete dimension distribution histogram in one go via the Hamming state tensor, namely by contracting $\set{M}^S$ with $\set{S}$ along all physical dimensions (Fig.~\ref{fig:contraction}). The array we seek results from the remaining free edge; see also Alg.~\ref{alg:dimension_distribution}.

\begin{algorithm}[ht]
\begin{algorithmic}[1]
\State Assemble $\set{M}^S$ with ranks $Q_0, Q_1, \dots , Q_{N-1}, Q_N = 1, N+1, \dots, N+1$ (Sec.~\ref{sec:mean_dimension})
\State Contract $\set{S}$ with $\set{M}^S$ as in Alg.~\ref{alg:mean_dimension}
\Return $\nu := $ resulting vector of $N$ elements
\end{algorithmic}
\caption{Given a Sobol TT $\set{S}$, compute the dimension distribution $\nu$.}
\label{alg:dimension_distribution}
\end{algorithm}

\subsection{Effective Dimension (Superposition Sense)} \label{sec:computing_effective_superposition}

The variance due to order $n$ and below is

\begin{equation}
\label{eq:overall_variance2}
\sum_{\pmb{\alpha} \given[\big] |\pmb{\alpha}| \le n} \set{S}_{\pmb{\alpha}} = \sum_{\pmb{\alpha}} (\set{S} \circ \set{M}^{\le n})_{\pmb{\alpha}} = <\set{S}, \set{M}^{\le n}>
\end{equation}

With Eq.~\ref{eq:overall_variance2} we can easily extract the superposed effective dimension: it is sufficient to iteratively find the smallest $k$ that yields a relative variance above the given threshold $1-\eps$; see also Alg.~\ref{alg:effective_superposition}.

\begin{algorithm}[ht]
\begin{algorithmic}[1]
\State Compute dimension distribution $\nu$ as in Alg.~\ref{alg:dimension_distribution}
\State $\bar{\nu} := \mathrm{cumulativeSum}(\nu)$
\For{$n = 1, \dots, N$}
	\If{$\bar{\nu}[n] \ge 1 - \eps$}
		\Return $d_S : = n$
	\EndIf
\EndFor
\end{algorithmic}
\caption{Given a Sobol TT $\set{S}$, compute its effective dimension (in the superposition sense) $d_S$ w.r.t. threshold $\eps$.}
\label{alg:effective_superposition}
\end{algorithm}

\subsection{Effective Dimension (Truncation Sense)} \label{sec:computing_effective_truncation}

The truncated effective dimension $d_T$ depends on the ordering of variables~\cite{GCAI:16}. However, the formula

\begin{equation}
\argmin_n \left\{ n \given[\big] \max_{\pmb{\alpha}} \left\{ (\set{S} \circ \set{M}^{\le n})_{\pmb{\alpha}} \right\} \ge 1-\eps \right\}
\end{equation}
gives us the truncated effective dimension w.r.t. the best possible ordering $\pmb{\alpha}$ of variables. The $n$ can be found iteratively (Alg.~\ref{alg:effective_truncation}).

\begin{algorithm}[ht]
\begin{algorithmic}[1]
\State Compute dimension distribution $\nu$ as in Alg.~\ref{alg:dimension_distribution}
\For{$n = 1, \dots, N$}
	\State Assemble $\set{M}^{\le n}$ (Sec.~\ref{sec:hamming_mask})
	\State $v := \max_{\pmb{\alpha}} \{\set{S} \circ \set{M}^{\le n}\}$ \Comment{Element-wise product using cross-approximation~\cite{OT:10}, tensor maximum found as in~\cite{MO:15}}
	\If{$v \ge 1 - \eps$}
		\Return $d_s : = n$
	\EndIf
\EndFor
\end{algorithmic}
\caption{Given a Sobol TT $\set{S}$, compute its effective dimension (in the truncation sense) $d_T$ w.r.t. threshold $\eps$ and best possible ordering of variables.}
\label{alg:effective_truncation}
\end{algorithm}

\subsection{Effective Dimension (Successive Sense)}

The following dot product is the total variance contributed by all tuples whose length is $n$ or less:

\begin{equation}
0 \le \left<\set{S}, \set{L}^{\le n}\right> \le 1
\end{equation}

\noindent
The effective dimension in the successive sense is therefore the smallest $n$ (Alg.~\ref{alg:effective_successive}):

\begin{equation}
d_s = \argmin_n \left\{ n \given[\big] \left<\set{S}, \set{L}^{\le n}\right> \ge 1-\eps \right\}
\end{equation}

\begin{algorithm}[ht]
\begin{algorithmic}[1]
\State $\vec{l} :=$ the contraction between $\set{S}$ and $\set{L}^S$ as in Alg.~\ref{alg:dimension_distribution}
\State $d_s := $ proceed as in Alg.~\ref{alg:effective_superposition} but using $\vec{l}$ instead of $\nu$
\Return $d_s$
\end{algorithmic}
\caption{Given a Sobol TT $\set{S}$, compute its effective dimension (in the successive sense) $d_s$ w.r.t. threshold $\eps$.}
\label{alg:effective_successive}
\end{algorithm}

\subsection{Shapley Values}

The $n$-th Shapley value can be interpreted as variance contributions in a SA context and computed in terms of the Sobol indices~\cite{Owen:14}. It is considered a challenging QoI from a computational point of view~\cite{Owen:14, Owen:17} since it relies on $2^{N-1}$ terms, each weighted by a combinatorial number. We observe that, in tensor form, it is equivalent to:

\begin{equation} \label{eq:computing_shapley}
\phi_n = \sum_{\pmb{\alpha} | n \in \pmb{\alpha}} \frac{S_{\pmb{\alpha}}}{|\pmb{\alpha}|} = (\set{S} \circ (1/\set{W}))^C_{\{n\}}
\end{equation}

This represents the closed version of the weighted Sobol indices, evaluated at tuples $\{1\}, ..., \{N\}$. Note that $1/\set{W}$ is a so-called \emph{Hilbert tensor}; such tensors are known to have high-rank but are extremely well compressible via low-rank expansions~\cite{OT:10}. We have confirmed this experimentally for this $2^N$ particular case (we first set $\set{W}_{\emptyset}$ to 1, to prevent division by 0): for $N = 50$, for example, a TT-rank of 7 is enough for $1/\set{W}$ to achieve a relative error under $0.0001\%$. In Alg.~\ref{alg:shapley_values} we show how to compute all Shapley values $\phi_1, \dots, \phi_N$ using Eq.~\ref{eq:computing_shapley}.

\begin{algorithm}[ht]
\begin{algorithmic}[1]
\State Compute $1/\set{W}$ \Comment{Using cross-approximation. Can be precomputed}
\State $\widehat{\set{S}} = [[\widehat{S}^{(1)}, \dots, \widehat{S}^{(N)}]] := \set{S} \circ (1/\set{W})$ \Comment{Denote its ranks as $1, Q_1, \dots, Q_{N-1}, 1$}
\For{$n = 1, \dots, N$} \Comment{Compute the closed tensor of a tensor~\cite{BPP:17}}
	\State $\widehat{\set{S}}^{{(n)}^C} := \mathrm{zeros}(Q_{n-1}, 2, Q_n)$
	\State $\widehat{\set{S}}^{{(n)}^C}[:, 0, :] := \widehat{\set{S}}^{(n)}[:, 0, :]$ \Comment{The first slice stays the same}
	\State $\widehat{\set{S}}^{{(n)}^C}[:, 1, :] := \widehat{\set{S}}^{(n)}[:, 0, :] + \widehat{\set{S}}^{(n)}[:, 1, :]$ \Comment{The second slice becomes the sum of both slices}
\EndFor
\State $\widehat{\set{S}}^C := [[\widehat{\set{S}}^{{(1)}^C}, \dots, \widehat{\set{S}}^{{(N)}^C}]]$
\For{$n = 1, \dots, N$}
	\State $\phi_n := \widehat{\set{S}}_{\{n\}}^C = \widehat{\set{S}}^C[\overbrace{0, \dots, 0}^{n-1}, 1, \overbrace{0, \dots, 0}^{N-n}]$
\EndFor
\Return $(\phi_1, \dots, \phi_N)$
\end{algorithmic}
\caption{Given a Sobol TT $\set{S}$, compute all $N$ Shapley values.}
\label{alg:shapley_values}
\end{algorithm}

\section{Experimental Results} \label{sec:results}

We have tested our metric computation algorithms\footnote{Our algorithms are provided within the Python package \url{https://github.com/rballester/ttrecipes} (the code for all models is available in the folder \texttt{examples/sensitivity\_analysis/}).} with three models of dimensionalities 10 and 20 on a desktop workstation (3.2GHz Intel i5). For each model, a TT surrogate with 100 points per axis is built using adaptive cross-approximation~\cite{OT:10, SO:11} as released in the Python library \emph{ttpy}~\cite{ttpy}, from which we extract its Sobol TT using the method from~\cite{BPP:17}. In all cases we compare all our resulting first-order Sobol and total Sobol indices with the quasi-MC algorithm by Saltelli et al.~\cite{SAACRT:10} implemented in the Sensitivity Analysis Library for Python (\emph{SALib})~\cite{HU:17}, and our Shapley values with the state-of-the-art MC method by Song et al.~\cite{SNS:16}. We also assess the accuracy of the proposed algorithm for the mean dimension $D_S$ by using the identity $D_S = \sum_n S^T_n$ due to Liu and Owen~\cite{LO:06}. We report the number of function evaluations for each method.
\subsection{Sobol ``G''  function}

Our first model is synthetic and widely used for benchmarking purposes in the SA literature:

\begin{equation}
\label{eq:sobol}
f(\vec{x}) := \prod_{n=1}^N \frac{|4 x_n - 2| + a_n}{1 + a_n}
\end{equation}
with $x_n \sim \mathcal{U}(0, 1)$ for all $n$. Coefficients $a_n$ are customizable and usually non-negative; we have chosen $a_1 := \dots := a_N := 0$ as in~\cite{BF:88} and~\cite{KFSM:11}. This model is partly provided as a sanity check: since $f$ is separable (i.e. has TT rank 1), we can expect to obtain an exact TT interpolator (up to axis discretization error). We took $N = 20$ dimensions. Due to variable symmetry, all Shapley values equal $1/N$. We furthermore know the analytical values for the SI and total SI from~\cite{KFSM:11}:

\begin{equation}
S_n = \frac{1}{3 \cdot \left(\left(\frac{4}{3}\right)^N - 1\right)}, \, S^T_n = \frac{S_n}{\left(\frac{4}{3}\right)^{1-N}} \mbox{ for } n = 1, \dots, N
\end{equation}

\begin{figure}[ht]\centering
	\includegraphics[width=0.7\columnwidth]{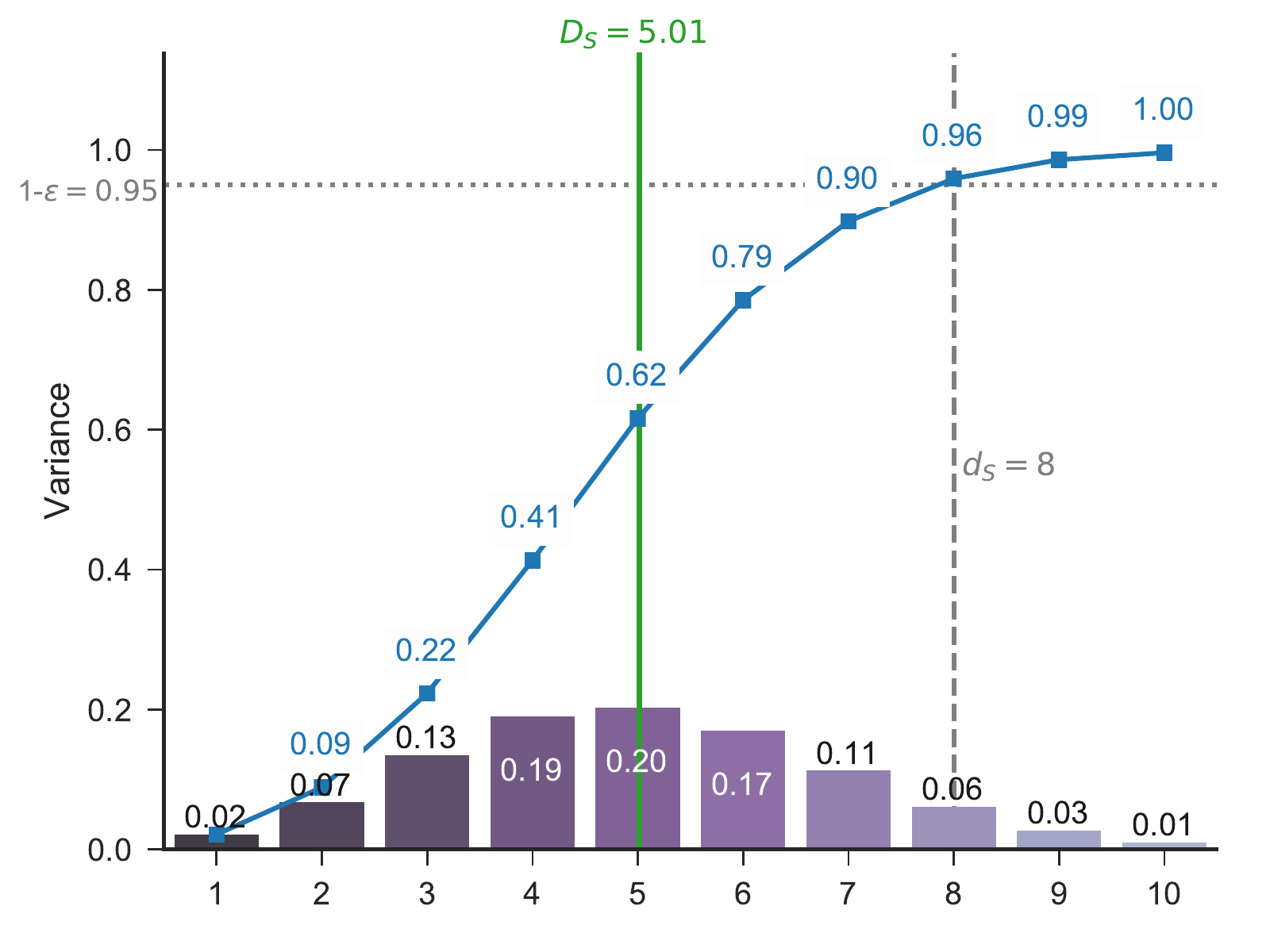}
	\caption {Sobol ``G'' function: dimension distribution (truncated after order $10$), mean dimension, and effective dimension in the superposition sense $d_S$ for $\eps = 0.05$.}
	\label{fig:sobolg_dim_dist_plot}
\end{figure}

Albeit separable, this model is high-dimensional and has high-order interactions, which makes it more challenging for MC- and QMC-based approaches. Our method used 59,400 function evaluations and required a total time of $3.8$s for building the surrogate model and computing all the sensitivity metrics. We run method~\cite{SAACRT:10} until its total SI relative error became $10\%$ of the ground-truth on average, for which $1.5 \times 10^8$ samples were needed. Method~\cite{SNS:16} for the Shapley values did not converge before $3 \times 10^8$ samples, after which the experiment was stopped. Numerical results are reported in Tab.~\ref{tab:sobolg1} (effective and mean dimensions), Tab.~\ref{tab:sobolg2} (Sobol/Shapley indices), and Fig.~\ref{fig:sobolg_dim_dist_plot} (dimension distribution). Note that the proposed method coincides with the ground-truth in all cases by at least three decimal digits. Note also that the QMC~\cite{SAACRT:10} only provides the first-order SI and total SI, whereas the Sobol TT contains indices for interactions of all possible order (although only 1st-order are shown).

\begin{table*}[ht]
	\centering
	\caption{Sobol ``G'' function: effective and mean dimensions.}
	
\begin{tabular}{ll}
\toprule
Dimension metric & Value \\
\midrule
Effective dimension ($\eps=0.05$) \\
\phantom{abc}\emph{Superposition} sense ($d_S$) & 8 (rel. var.: $0.959$) \\
\phantom{abc}\emph{Truncation} sense ($d_T$) $[$ $X_{1}$ $\ \cdots \ $ $X_{20}$ $]$ & 20 (rel. var.: $1.000$) \\
\phantom{abc}\emph{Successive} sense ($d_s$) & 20 (rel. var.: $1.000$) \\
\midrule
Mean dimension ($D_s$) &  \\
\phantom{abc}Exact & 5.016 \\
\phantom{{abc}}[Ours] Automata & 5.014 \\
\phantom{{abc}}[Ours] $\sum_n S^T_n$ & 5.014 \\
\phantom{{abc}}[QMC~\cite{SAACRT:10}] $\sum_n S^T_n$ & 4.720 \\

\bottomrule
\end{tabular}

	\label{tab:sobolg1}
\end{table*}

\begin{table*}[ht]
	\centering
	\caption{Sobol ``G'' function: Shapley values and Sobol indices of all inputs. The Monte Carlo (MC) method~\cite{SNS:16} was stopped without having converged at $3 \times 10^8$ samples. The quasi-Monte Carlo (QMC) algorithm~\cite{SAACRT:10} used $1.5 \times 10^8$ samples.}
	\resizebox{1\columnwidth}{!}{
		
\begin{tabular}{lccccccccc}
\toprule
  \multirow{2}{*}{ } & \multicolumn{3}{c}{Shapley value} & \multicolumn{3}{c}{Sobol index} & \multicolumn{3}{c}{Total Sobol index} \\
\cmidrule(r){2-4} \cmidrule(r){5-7}  \cmidrule(r){8-10}
   & Ours & Exact & MC~\cite{SNS:16} & Ours & Exact & QMC~\cite{SAACRT:10} & Ours & Exact & QMC~\cite{SAACRT:10}\\ 
\midrule

$X_{1}$ & 0.050 & 0.050 & - & 0.001 & 0.001 & 0.001 & 0.251 & 0.251 & 0.257 \\
$X_{2}$ & 0.050 & 0.050 & - & 0.001 & 0.001 & 0.001 & 0.251 & 0.251 & 0.226 \\
$X_{3}$ & 0.050 & 0.050 & - & 0.001 & 0.001 & 0.001 & 0.251 & 0.251 & 0.198 \\
$X_{4}$ & 0.050 & 0.050 & - & 0.001 & 0.001 & 0.001 & 0.251 & 0.251 & 0.238 \\
$X_{5}$ & 0.050 & 0.050 & - & 0.001 & 0.001 & 0.001 & 0.251 & 0.251 & 0.277 \\
$X_{6}$ & 0.050 & 0.050 & - & 0.001 & 0.001 & 0.001 & 0.251 & 0.251 & 0.222 \\
$X_{7}$ & 0.050 & 0.050 & - & 0.001 & 0.001 & 0.001 & 0.251 & 0.251 & 0.209 \\
$X_{8}$ & 0.050 & 0.050 & - & 0.001 & 0.001 & 0.001 & 0.251 & 0.251 & 0.237 \\
$X_{9}$ & 0.050 & 0.050 & - & 0.001 & 0.001 & 0.001 & 0.251 & 0.251 & 0.272 \\
$X_{10}$ & 0.050 & 0.050 & - & 0.001 & 0.001 & 0.001 & 0.251 & 0.251 & 0.231 \\
$X_{11}$ & 0.050 & 0.050 & - & 0.001 & 0.001 & 0.001 & 0.251 & 0.251 & 0.210 \\
$X_{12}$ & 0.050 & 0.050 & - & 0.001 & 0.001 & 0.001 & 0.251 & 0.251 & 0.184 \\
$X_{13}$ & 0.050 & 0.050 & - & 0.001 & 0.001 & 0.001 & 0.251 & 0.251 & 0.263 \\
$X_{14}$ & 0.050 & 0.050 & - & 0.001 & 0.001 & 0.002 & 0.251 & 0.251 & 0.217 \\
$X_{15}$ & 0.050 & 0.050 & - & 0.001 & 0.001 & 0.001 & 0.251 & 0.251 & 0.279 \\
$X_{16}$ & 0.050 & 0.050 & - & 0.001 & 0.001 & 0.001 & 0.251 & 0.251 & 0.239 \\
$X_{17}$ & 0.050 & 0.050 & - & 0.001 & 0.001 & 0.001 & 0.251 & 0.251 & 0.240 \\
$X_{18}$ & 0.050 & 0.050 & - & 0.001 & 0.001 & 0.001 & 0.251 & 0.251 & 0.249 \\
$X_{19}$ & 0.050 & 0.050 & - & 0.001 & 0.001 & 0.001 & 0.251 & 0.251 & 0.217 \\
$X_{20}$ & 0.050 & 0.050 & - & 0.001 & 0.001 & 0.001 & 0.251 & 0.251 & 0.254 \\

\bottomrule
\end{tabular}

	}
	\label{tab:sobolg2}
\end{table*}

\subsection{Simulated decay chain}

The second example simulates a radioactive decay chain that concatenates Poisson processes for 11 chemical species. It is a linear \emph{Jackson network}, i.e. each species (except the last one) can decay into the next species in the chain. We model 10 parameters, namely the decay rates $\lambda_n$ of the 10 first species. The simulation result $f_T(\lambda_1, \dots, \lambda_{10})$ is the amount of stable material (last node in the chain) measured after a certain time span $T$. To evaluate each sample of $f_T$ we simulate the chain by discretizing the span $T$ into timesteps of one day. The $\lambda_n$ represent the fraction of each material that decays every day. They are independent and uniformly distributed in the interval $[0.00063281, 0.00756736]$, which corresponds to half-lives from 3 years down to 3 months. 


The resulting effective and mean dimension metrics for a time span of $T = 2$ years are reported in Tab.~\ref{tab:decay1}, and the Sobol and Shapley sensitivity indices in Tab.~\ref{tab:decay2}. Our method used 134,400 function evaluations and took a total time of $8.4$s. Methods~\cite{SNS:16} and~\cite{SAACRT:10} were run with $1.4 \times 10^7$ and 600,000 samples respectively. These were the minimal numbers such that their suggested confidence intervals (one standard error for~\cite{SAACRT:10}, two for~\cite{SNS:16}) were on average the $10\%$ of their estimated indices. Note that our results converge unanimously to $S_1 = \dots = S_{10}$, $S^T_1 = \dots = S^T_{10}$, and $\phi_1 = \dots = \phi_{10} = 1/10$.

\begin{table*}[ht]
	\centering
	\caption{Decay chain: effective and mean dimensions (time span = 2 years)}
	
\begin{tabular}{ll}
\toprule
Dimension metric & Value \\
\midrule
Effective dimension ($\eps=0.05$) \\
\phantom{abc}\emph{Superposition} sense ($d_S$) & 3 (rel. var.: $0.959$) \\
\phantom{abc}\emph{Truncation} sense ($d_T$) $[$ $\lambda_{1}$ $\ \cdots \ $ $\lambda_{10}$ $]$ & 10 (rel. var.: $1.000$) \\
\phantom{abc}\emph{Successive} sense ($d_s$) & 9 (rel. var.: $0.978$) \\
\midrule
Mean dimension ($D_s$) &  \\
\phantom{{abc}}[Ours] Automata & 1.728 \\
\phantom{{abc}}[Ours] $\sum_n S^T_n$ & 1.728 \\
\phantom{{abc}}[QMC~\cite{SAACRT:10}] $\sum_n S^T_n$ & 1.722 \\

\bottomrule
\end{tabular}

	\label{tab:decay1}
\end{table*}

\begin{table*}[ht]
	\centering
	\caption{Decay chain: Shapley values and Sobol indices (time span = 2 years)}
	
\begin{tabular}{lcccccc}
\toprule
  \multirow{2}{*}{Variable} & \multicolumn{2}{c}{Shapley value} & \multicolumn{2}{c}{Sobol index} & \multicolumn{2}{c}{Total Sobol index} \\
\cmidrule(r){2-3} \cmidrule(r){4-5}  \cmidrule(r){6-7}
   & Ours & MC~\cite{SNS:16} & Ours & QMC~\cite{SAACRT:10} & Ours & QMC~\cite{SAACRT:10} \\
\midrule

$\lambda_{1}$ & 0.100 & 0.096 & 0.049 & 0.049 & 0.173 & 0.171 \\
$\lambda_{2}$ & 0.100 & 0.099 & 0.049 & 0.052 & 0.173 & 0.171 \\
$\lambda_{3}$ & 0.100 & 0.093 & 0.049 & 0.050 & 0.173 & 0.170 \\
$\lambda_{4}$ & 0.100 & 0.102 & 0.049 & 0.049 & 0.173 & 0.171 \\
$\lambda_{5}$ & 0.100 & 0.100 & 0.049 & 0.048 & 0.173 & 0.167 \\
$\lambda_{6}$ & 0.100 & 0.099 & 0.049 & 0.050 & 0.173 & 0.184 \\
$\lambda_{7}$ & 0.100 & 0.111 & 0.049 & 0.051 & 0.173 & 0.171 \\
$\lambda_{8}$ & 0.100 & 0.109 & 0.049 & 0.051 & 0.173 & 0.170 \\
$\lambda_{9}$ & 0.100 & 0.096 & 0.049 & 0.048 & 0.173 & 0.172 \\
$\lambda_{10}$ & 0.100 & 0.097 & 0.049 & 0.050 & 0.173 & 0.176 \\

\bottomrule
\end{tabular}

	\label{tab:decay2}
\end{table*}

\begin{figure}[ht]\centering
	\includegraphics[width=\columnwidth]{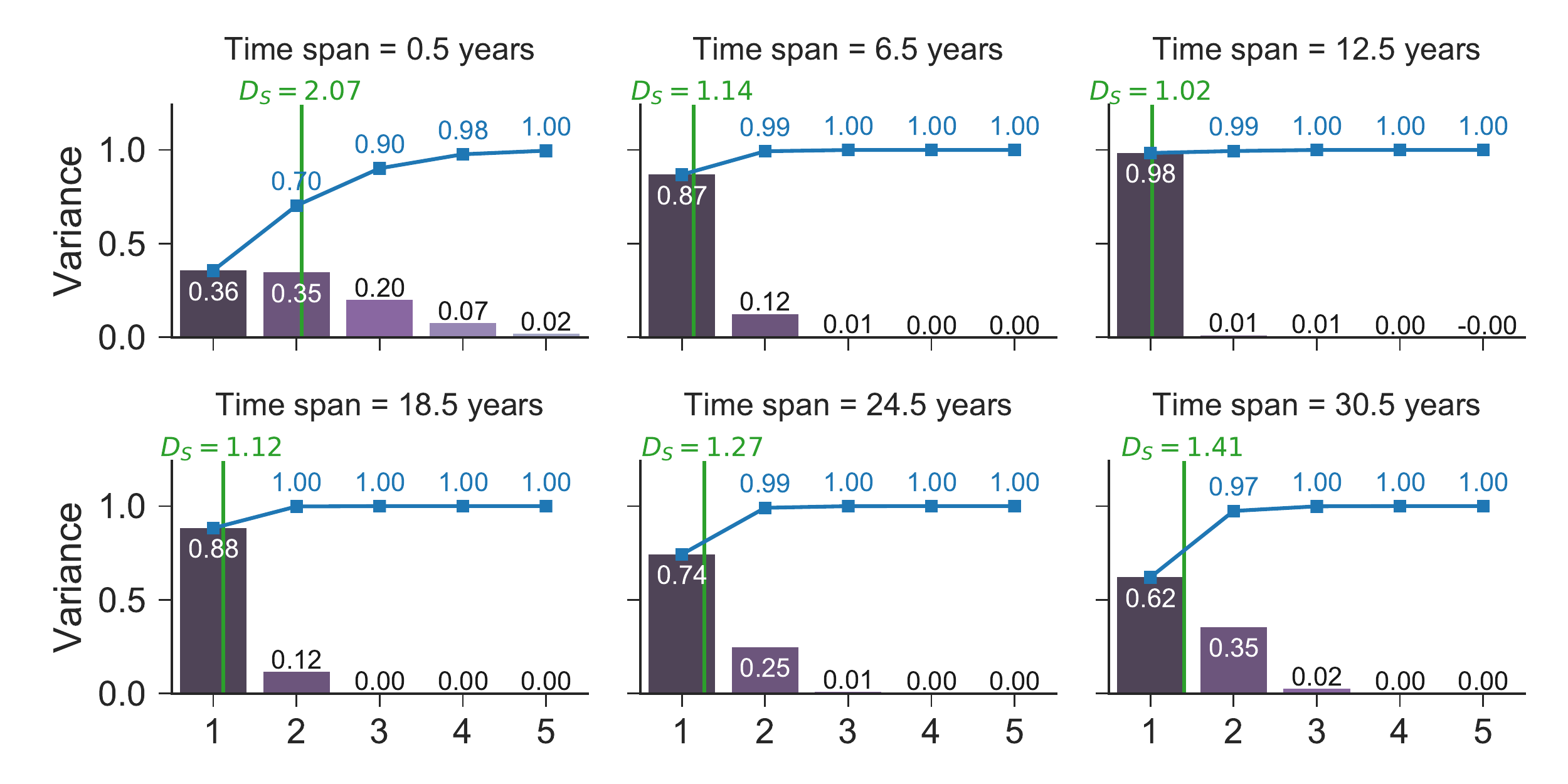}
	\caption{Decay chain: mean dimension and dimension distribution (truncated after order $5$) for 6 different time spans $0.5 \le T \le 30$. Very short and very long time spans result in higher-order interactions.}
	\label{fig:decay_series_dim_dist_plot}
\end{figure}

\begin{figure}[ht]\centering
	\includegraphics[width=0.9\columnwidth]{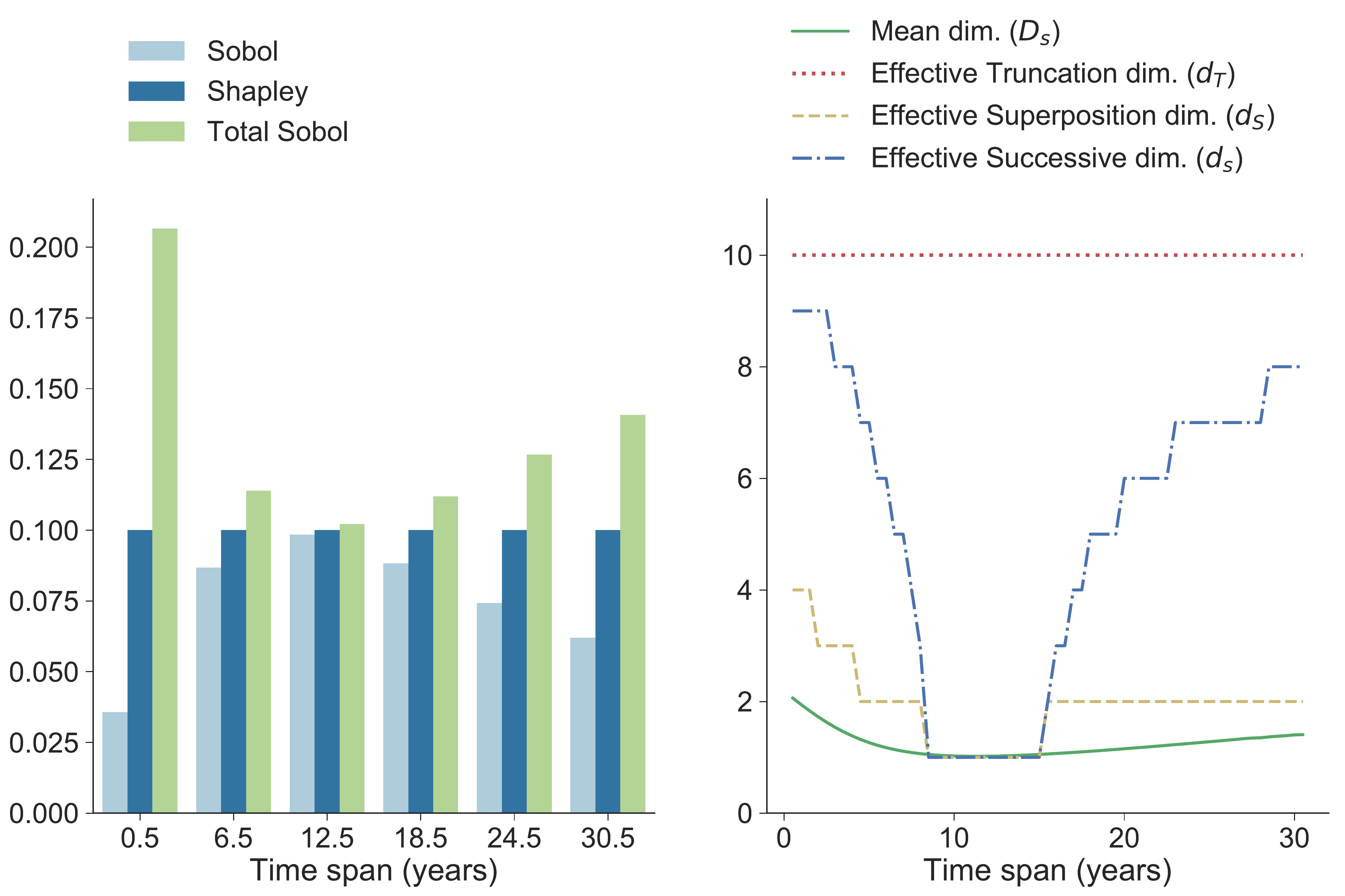}
	\caption{Decay chain: all metrics as a function of the time span.
	Note (left plot) how the SI and total SI follow an inverse relationship with one another and meet at about $T = 12$ years, while the Shapley values are constant.}
	\label{fig:decay_series_mixed}
\end{figure}

We have also studied the sensitivity behavior when varying the simulated time span: $0.5 \le T \le 30.5$ years (see dimension distributions in Fig.~\ref{fig:decay_series_dim_dist_plot}). Fig.~\ref{fig:decay_series_mixed} depicts the evolution of the Shapley values, Sobol and total indices and the dimension metrics within the same range. Note how the average interaction order is higher for extreme values of $T$ and lower in the center ($\approx 1$ at around 12 years). This behavior hints that, for very short or long time spans, several decay rates need to have specific values to affect the output of the function. In the former case, for example, no amount of the final species is detected unless many decay rates are fast enough. On the other hand, after a very long time span $T$, all materials have decayed completely unless several rates are slow enough. The successive dimension $d_s$ (right chart in Fig.~\ref{fig:decay_series_mixed}) strongly reflects this behavior as well: the metric is highly sensitive to such high-order consecutive interactions, and is therefore useful for this kind of time-series models.

Finally, this example illustrates how the Shapley values $\phi_n$ are able to recognize the equal importance of all variables and act as the best overall summary of their influence. This is consistent with Owen~\cite{Owen:14}, their original proponent in the context of sensitivity analysis.

\subsection{Fire-spread Model}

Last we model the rate of fire spread in the Mediterranean shrublands according to 10 variables that are fed into Rothermel's equations~\cite{Rothermel:72}. 
Both~\cite{SPTP:01} and~\cite{SNS:16} have analyzed this model from the sensitivity analysis perspective. We compare our results with the latter work for the case of independently (but non-uniformly) distributed parameters. We use the updated equations as in~\cite{SNS:16}, which incorporates the modifications by Albini~\cite{Albini:76} and Catchpole and Catchpole~\cite{CC:91}, with the marginal PDFs shown in Tab.~\ref{tab:fire_params}.

\begin{figure}[ht]\centering
	\includegraphics[width=0.5\columnwidth]{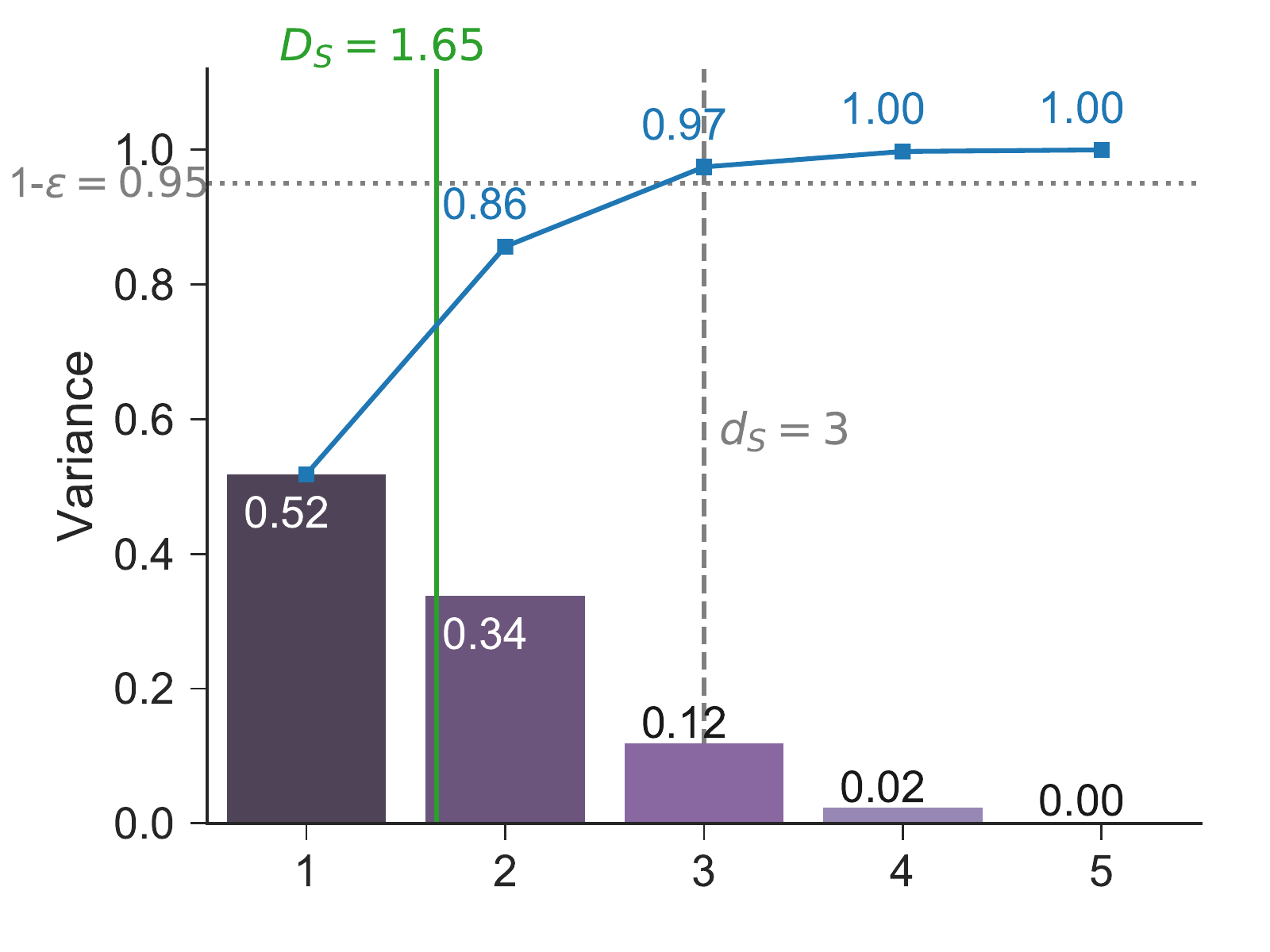}
	\caption {Fire spread: dimension distribution (truncated after order $5$).}
	\label{fig:fire_dim_dist_plot}
\end{figure}

Our method used $4.58 \times 10^6$ samples and took $9.2$s to extract all sensitivity metrics. Method~\cite{SNS:16} was run with $4.6 \times10^7$ samples (just as in the original paper), while we again run method~\cite{SAACRT:10} so as to attain a $10\%$ confidence interval ($6.96 \times 10^6$ samples). Results are reported in Tab.~\ref{tab:fire1} (effective and mean dimensions) and Tab.~\ref{tab:fire2} (Shapley values and Sobol indices) as well as Fig.~\ref{fig:fire_dim_dist_plot} (dimension distribution). The analytical metrics for this model are unknown, but our results are in all cases within the other two methods' intervals of confidence.

\begin{table*}[ht]
	\centering
	\caption{Fire spread: parameter description and marginal PDFs~\cite{SNS:16}.}
\resizebox{1\columnwidth}{!}{
		\begin{tabular}{cll}
		\toprule
		Variable & Description  & Distribution \\
		\midrule
		$\delta$ & Fuel depth ($\mathrm{cm}$) & $Log\mathcal{N}(2.19, 0.517) \cap [5, \infty)$ \\
		$\sigma$ & Fuel particle area-to-volume ratio ($1/\mathrm{cm}$) & $Log\mathcal{N}(3.31, 0.294)$ \\
		$h$ & Fuel particle low heat content ($\mathrm{Kcal/kg}$) & $Log\mathcal{N}(8.48, 0.063)$ \\
		$\rho_p$ & Oven-dry particle density ($\mathrm{D.W.g/cm^3}$) & $Log\mathcal{N}(-0.592, 0.219)$ \\
		$m_l$ & Moisture content of the live fuel ($\mathrm{H_{2}O\ g / D.W.g}$) & $\mathcal{N}(1.18, 0.377) \cap [0, \infty)$ \\
		$m_d$ & Moisture content of the dead fuel ($\mathrm{H_{2}O\ g / D.W.g}$) & $\mathcal{N}(0.19, 0.047)$ \\
		$S_T$ & Fuel particle total mineral content ($\mathrm{MIN.g / D.W.g}$) & $\mathcal{N}(0.049, 0.011) \cap [0, \infty)$ \\
		$U$ & Wind speed at midflame height ($\mathrm{km/h}$) & $6.9 \cdot Log\mathcal{N}(1.0174, 0.5569)$ \\
		$\tan \phi$ & Slope  $ $ & $\mathcal{N}(0.38, 0.186) \cap [0, \infty)$ \\
		$P$ & Dead to total fuel loading ratio $ $ & $Log\mathcal{N}(-2.19, 0.64) \cap (-\infty, 1]$ \\
		\bottomrule
	\end{tabular}
}
	\label{tab:fire_params}
\end{table*}

\begin{table*}[ht]
	\centering
	\caption{Fire spread: effective and mean dimensions.}
	
\begin{tabular}{ll}
\toprule
Dimension metric & Value \\
\midrule
Effective dimension ($\eps=0.05$) \\
\phantom{abc}\emph{Superposition} sense ($d_S$) & 3 (rel. var.: $0.975$) \\
\phantom{abc}\emph{Truncation} sense ($d_T$) $[$ $\delta$, $\sigma$, $m_l$, $m_d$, $U$, $P$ $]$ & 6 (rel. var.: $0.961$) \\
\phantom{abc}\emph{Successive} sense ($d_s$) & 8 (rel. var.: $0.963$) \\
\midrule
Mean dimension ($D_s$) &  \\
\phantom{{abc}}[Ours] Automata & 1.653 \\
\phantom{{abc}}[Ours] $\sum_n S^T_n$ & 1.653 \\
\phantom{{abc}}[QMC~\cite{SAACRT:10}] $\sum_n S^T_n$ & 1.714 \\

\bottomrule
\end{tabular}

	\label{tab:fire1}
\end{table*}

\begin{table*}[ht]
	\centering
	\caption{Fire spread: Shapley values and Sobol indices.}
	
\begin{tabular}{lcccccc}
\toprule
  \multirow{2}{*}{Variable} & \multicolumn{2}{c}{Shapley value} & \multicolumn{2}{c}{Sobol index} & \multicolumn{2}{c}{Total Sobol index} \\
\cmidrule(r){2-3} \cmidrule(r){4-5}  \cmidrule(r){6-7}
   & Ours & MC~\cite{SNS:16} & Ours & QMC~\cite{SAACRT:10} & Ours & QMC~\cite{SAACRT:10} \\
\midrule

$\delta$ & 0.203 & 0.217 & 0.106 & 0.106 & 0.331 & 0.348 \\
$\sigma$ & 0.125 & 0.132 & 0.048 & 0.048 & 0.236 & 0.228 \\
$h$ & 0.003 & -0.020 & 0.001 & 0.001 & 0.005 & 0.006 \\
$\rho_p$ & 0.012 & 0.013 & 0.004 & 0.005 & 0.023 & 0.023 \\
$m_l$ & 0.231 & 0.231 & 0.142 & 0.143 & 0.347 & 0.364 \\
$m_d$ & 0.165 & 0.183 & 0.095 & 0.097 & 0.259 & 0.262 \\
$S_T$ & 0.002 & -0.006 & 0.001 & 0.001 & 0.005 & 0.004 \\
$U$ & 0.202 & 0.210 & 0.090 & 0.093 & 0.354 & 0.387 \\
$\tan \phi$ & 0.004 & -0.014 & 0.002 & 0.002 & 0.006 & 0.006 \\
$P$ & 0.053 & 0.054 & 0.029 & 0.029 & 0.087 & 0.086 \\

\bottomrule
\end{tabular}

	\label{tab:fire2}
\end{table*}

\section{Discussion and Conclusions} \label{sec:conclusions}

The Sobol indices have long been in the spotlight of the variance-based sensitivity analysis and uncertainty quantification communities. In recent years, several other sensitivity metrics are becoming popular that capitalize on these building blocks and are combinatorial in nature. State-of-the-art algorithms and sampling schemes that estimate these advanced metrics often require rather large numbers of samples, i.e. in the order of millions of function evaluations. Therefore, modern methods are growingly dependent on fast analytical or surrogate models in order to keep overall computation times low. We argue that similar or lower sampling budgets are actually sufficient for building highly accurate TT surrogates. Once such a tensor representation is available, its advantages become evident: a wide range of operations can be computed expeditiously in the compressed domain, including differentiation, integration, optimization, element-wise dot products and functions, and in particular fast evaluation of many sensitivity metrics as proposed in this paper. Our simulations have demonstrated two aspects of the proposed methods:

\begin{itemize}
	\item Correctness: several families of models admit an exact low-rank TT decomposition, and many others can be well approximated by one. In addition, all proposed TT automata are error-free. Once we have a high-quality surrogate, we can extract metrics that are close to analytical metrics by several decimal digits of precision.
	\item Efficiency: cross-approximation is efficient in the sense that it needs a number of samples (model evaluations) that is proportional to the decomposition's compressed size. When the underlying model can be well approximated by a low-rank TT tensor, this leads to relatively low sampling budget requirements. The usual advantages of surrogate modeling apply: once we have our TT tensor (and especially its Sobol TT), arbitrary metrics can be computed from it on demand without further sampling.
\end{itemize}

To summarize, we have contributed algorithms for advanced variance-based sensitivity analysis that leverage tensor decompositions heavily, in particular the TT model. Such decompositions can represent exponential numbers of indices and QoIs in an extremely compact manner. Throughout our paper we exploited three crucial tools. First, the Sobol TT is able to gather all Sobol indices. While most surrogate-based SA approaches use a metamodel as a departure point for metric computation, we go one step further and work with a data structure that already contains all elementary indices of interest. Second, the automaton-inspired TT tensors contributed here are in charge of selecting and weighing index tuples as needed, and they all have a moderate rank. Last, tensor-tensor contractions (i.e. dot products) in the TT format have a polynomial cost w.r.t. number of dimensions and, combined with the previous elements, can produce sophisticated QoIs in a matter of few sequential steps. Other auxiliary (but still important) tools include adaptive cross-approximation, which is useful for building TT surrogates and auxiliary tensors, and global optimization in the compressed domain.

\subsection*{Limitations and Future Work}

As we have just highlighted and demonstrated in this work, moderate to large sampling budgets in black-box settings are usually sufficient for building high-quality TT models via ACA. However, the particular case remains of how to extract reliable sensitivity metrics when limited by more modest sampling budgets, e.g. of a few hundreds or thousands of samples. Using an intermediate surrogate often becomes the obligated approach in these cases, also for the other state-of-the-art methods~\cite{OCMB:15}. On the TT side, one may either run ACA on the intermediate model or build directly the tensor from a limited set of samples using e.g. low-rank tensor completion techniques. Research on TT completion usually strives for a generalization error as low as possible, e.g. as estimated by cross-validation strategies. For SA we must additionally place an emphasis on consistence, i.e. on how to determine and enforce the conditions under which TT surrogates trained on the same limited ground-truth samples, but using different methods, yield the same sensitivity metrics (as it is desirable). These aspects will be the subject of future research.


\section*{Acknowledgments}

This work was partially supported by the University of Zurich's Forschungskredit ``Candoc'' (grant number FK-16-012). The authors wish to thank Eunhye Song for providing code for the fire-spread model~\cite{SNS:16}.

\bibliographystyle{siamplain}
\bibliography{references}
\end{document}